\documentclass{proc-l}

\newtheorem{theorem}{Theorem}[section]

\theoremstyle{definition}

\newtheorem{property}[theorem]{Property}

\theoremstyle{remark}

\def\C{\mathbb{C}}
\def\N{\mathbb{N}}

\def\Li{\textrm{Li}_2}

\numberwithin{equation}{section}
\DeclareRobustCommand{\qbinom}{\genfrac[]{0pt}{}}
\begin{document}

\title[Euler's $q$-logarithm]{Leonhard Euler and a \boldmath$q$-analogue of the logarithm}

\dedicatory{at the 300th anniversary of Euler's birth}

\author{Erik Koelink}
\address{IMAPP, FNWI, Radboud Universiteit, Toernooiveld 1, 6525 ED Nijmegen, the Netherlands}
\email{e.koelink@math.ru.nl}

\author{Walter Van Assche}
\address{Departement Wiskunde, 
Katholieke Universiteit Leuven,
Celestijnenlaan 200B,
B-3001 Leuven, Belgium}
\email{walter@wis.kuleuven.be}
\thanks{WVA supported by research grant OT/04/21 of Katholieke Universiteit Leuven, research project G.0455.04 of FWO-Vlaanderen, and INTAS research network 03-51-6637}


\date{\today}

\commby{Peter Clarkson}

\begin{abstract}
We study a $q$-logarithm which was introduced by Euler and give
some of its properties. This $q$-logarithm did not get much attention
in the recent literature. We derive basic properties, some of which 
were already given by Euler in a 1751-paper and 1734-letter to Daniel
Bernoulli. The corresponding $q$-analogue of the dilogarithm is 
introduced. The relation to the values at $1$ and $2$ of a $q$-analogue
of the zeta function is given. 
We briefly describe some other $q$-logarithms
that have appeared in the recent literature. 
\end{abstract}

\maketitle

\section{Introduction}
In a paper from 1751, Leonhard Euler (1707--1783) introduced the series \cite[\S 6]{euler}
\begin{equation}  \label{eq:s}
     s = \sum_{k=1}^\infty  \frac{(1-x)(1-x/a)\cdots(1-x/a^{k-1})}{1-a^k}.
\end{equation} 
We will take $q=1/a$, 
then this series is convergent for $|q| <1$ and $x \in \mathbb{C}$.
 In this paper we will assume $0 < q < 1$.
Then this becomes
\begin{equation}  \label{eq:S}
    S_q(x) = - \sum_{k=1}^\infty  \frac{q^k}{1-q^k} (x;q)_k, 
\end{equation}
where $(x;q)_0=1$, $(x;q)_k = (1-x)(1-xq)\cdots(1-xq^{k-1})$. 
This can be written as a basic hypergeometric series
\[   S_q(x) = - \frac{q(1-x)}{1-q}\ {}_3\phi_2 \left( \genfrac{.}{.}{0pt}{}{q,q,qx}
{q^2,0}\ ;q,q \right). \]
Euler had come across this series much earlier in an attempt to
interpolate the logarithm at powers $a^k$ (or $q^{-k}$), see, e.g., Gautschi's
comment \cite{gautschi} discussing Euler's letter 
to Daniel Bernoulli where Euler introduced the function for $a=10$. 
Euler was aware that this interpolation did
not work very well, see \cite[\S 3-4]{gautschi}.
The function in \eqref{eq:S} does not seem to appear
in the recent literature, even though it has some nice properties. We will prove some
of its properties, some already obtained by Euler \cite{euler}, 
and indicate why this should be called a $q$-analogue of the logarithm.
A first reason is that for $0 < q < 1$
\[  \lim_{q \to 1} (1-q) S_q(x) = -\sum_{k=1}^\infty \lim_{q \to 1} q^k \frac{1-q}{1-q^k} (x;q)_k =
   - \sum_{k=1}^\infty \frac{(1-x)^k}{k} = \log x, \]
which is only a formal limit transition, since interchanging limit and sum
seems hard to justify. 

In Sections \ref{sec:Sentire}--\ref{sec:qdiffeq}
we study this $q$-analogue of the logarithm more closely.
In particular, we reprove some of Euler's results. Then we go on to
extend the definition in Section \ref{sec:extensionqlog}. Finally, we study the corresponding
$q$-analogue of the dilogarithm in Section 
\ref{sec:qdilog}. It involves also the values at $1$ and $2$
of a $q$-analogue of the $\zeta$-function. 
We give a (incomplete) list of some other
$q$-analogues of the logarithm appearing in the literature in
Section \ref{sec:otherqlogs}. 
The purpose of this note is to draw attention to the $q$-analogues of the logarithm, dilogarithm and
$\zeta$-function for which we expect many interesting results remain to be discovered.

Many results in this note use the $q$-binomial theorem \cite[\S 1.3]{GR}, \cite[\S 10.2]{AAR}
\begin{equation}  \label{eq:qbin}
 \frac{(ax;q)_\infty}{(x;q)_\infty} = \sum_{j=0}^\infty \frac{(a;q)_j}{(q;q)_j} x^j,
 \qquad |x| < 1.  
\end{equation}
We also use the $q$-exponential functions \cite[p.~9]{GR}, \cite[p.~492]{AAR}
\begin{eqnarray*}
   && e_q(z) = \frac{1}{(z;q)_\infty} = \sum_{n=0}^\infty \frac{z^n}{(q;q)_n}, \qquad |z| < 1, \\
   && E_q(z) = (-z;q)_\infty = \sum_{n=0}^\infty \frac{q^{n(n-1)/2}}{(q;q)_n} \ z^n .  
\end{eqnarray*}

\textit{Acknowledgement} This paper was triggered by a lecture
on the 1734-letter of Leonhard Euler to Daniel Bernoulli by Walter Gautschi, see
\cite{gautschi}. We thank Walter Gautschi for useful discussions and for providing
us with a translation of \cite{euler} (the 
translation can be downloaded from the E190 page of the Euler archive at 
\texttt{http://www.math.dartmouth.edu/\~{\,}euler}).  We thank the referee for
useful comments. The work of EK for this paper was mainly done at Technische
Universiteit Delft.

\section{The $q$-logarithm as an entire function}\label{sec:Sentire}

First of all we will show that the function $S_q$ in (\ref{eq:S}) is an entire
function, and as such it is a nicer function than the logarithm, which has a cut
along the negative real axis.

\begin{property}\label{prop:entireS}
The function $S_q$ defined in (\ref{eq:S}) is an entire function of order zero.
\end{property}

\begin{proof}
For $k \in \mathbb{N}$ the $q$-Pochhammer $(z;q)_k$ is a polynomial of degree $k$
with zeros at $1,1/q,\ldots,1/q^{k-1}$. For $|z| \leq r$ we have the simple bound
\[   |(z;q)_k| \leq (1+r)(1+r|q|)\cdots(1+r|q|^{k-1}) = (-r;|q|)_k < (-r;|q|)_\infty  \]
and hence the partial sums are uniformly bounded on the ball $|z| \leq r$:
\[  \left| - \sum_{k=1}^n \frac{q^k}{1-q^k} (z;q)_k \right|  \leq (-r;|q|)_\infty  \sum_{k=1}^\infty \frac{|q|^k}{1-|q|^k}.
\]
The partial sums therefore are a normal family and are uniformly convergent 
on every compact subset of the complex plane. 
The limit of these partial sums is $S_q(z)$ and is therefore an entire function of the complex
variable $z$.

Let $M(r)=\max_{|z|\leq r} |S_q(z)|$, then 
\[   M(r) \leq (-r;|q|)_\infty  \sum_{k=1}^\infty \frac{|q|^k}{1-|q|^k}  \]
and $(-r;|q|)_\infty = E_{|q|}(r)$ is the maximum of $E_{|q|}(z)$ on the ball $\{|z| \leq r\}$. 
The function $E_q$ is an entire function of order zero, which can be seen from the coefficients
$a_n$ of its Taylor series and the formula \cite[Theorem 2.2.2]{boas}
\begin{equation}  \label{eq:boas}
   \limsup_{n\to \infty} \frac{n\log n}{\log(1/|a_n|)}  
\end{equation}
for the order of $\sum_{n=0}^\infty a_n z^n$. Hence also $S_q$ has order zero. 
\end{proof}

Observe that for $0 < q < 1$ we have
\[     M(r) = \max_{|z|\leq r} |S_q(z)| = \sum_{k=1}^\infty \frac{q^k}{1-q^k} (-r;q)_k  \]
and some simple bounds give
\[   (q;q)_\infty   \sum_{k=1}^\infty \frac{q^k}{(q;q)_k} (-r;q)_k
                        \leq  M(r) \leq (-r;q)_\infty \sum_{k=1}^\infty \frac{q^k}{1-q^k} . \]
For the lower bound we can use the $q$-binomial theorem (\ref{eq:qbin}) to find
\[     (-rq;q)_\infty- (q;q)_\infty  \leq M(r) \leq (-r;q)_\infty \sum_{k=1}^\infty \frac{q^k}{1-q^k} \]
which shows that $M(r)$ behaves like $E_q(qr)-C_1 \leq M(r) \leq C_2 E_q(r)$, where $C_1$ and $C_2$ are
constants (which depend on $q$).

Euler \cite[\S 14-15]{euler} essentially also stated the following Taylor expansion. 

\begin{property}\label{prop:TaylorS}
The $q$-logarithm (\ref{eq:S}) has the following Taylor series around $x=0$:
\[    S_q(x) = - \sum_{k=1}^\infty \frac{q^k}{1-q^k} \left( 1 + q^{k(k-1)/2} \frac{(-x)^k}{(q;q)_k}
   \right).  \]
\end{property}

\begin{proof} Use the $q$-binomial theorem (\ref{eq:qbin}) with $x=zq^k$ and $a=q^{-k}$ to find
\begin{equation}  \label{qbinfin}
    (z;q)_k = \sum_{j=0}^k \qbinom{k}{j} q^{j(j-1)/2} (-z)^j, 
\quad \qbinom{k}{j} = \frac{(q;q)_k}{(q;q)_j(q;q)_{k-j}}.  
\end{equation}
Use this in (\ref{eq:S}), and change the order of summation to find
\[    S_q(x) = - \sum_{k=1}^\infty \frac{q^k}{1-q^k} - \sum_{j=1}^\infty q^{j(j-1)/2} (-x)^j
     \sum_{k=j}^\infty \frac{q^k}{1-q^k} \frac{(q;q)_k}{(q;q)_j (q;q)_{k-j}}.  \]
With a new summation index $k=j+\ell$ this becomes
\[  S_q(x) = - \sum_{k=1}^\infty \frac{q^k}{1-q^k} - \sum_{j=1}^\infty \frac{q^j}{1-q^j}
     q^{j(j-1)/2} (-x)^j \sum_{\ell=0}^\infty q^\ell \frac{(q^j;q)_{\ell}}{(q;q)_{\ell}}.  \]
Now use the $q$-binomial theorem \eqref{eq:qbin} to sum over $\ell$ to find 
\[   S_q(x) = - \sum_{k=1}^\infty \frac{q^k}{1-q^k} - \sum_{j=1}^\infty \frac{q^j}{1-q^j}
      q^{j(j-1)/2}  \frac{(-x)^j}{(q;q)_j}, \]
and if we combine both series then the required expansion follows.  
\end{proof}

This result can be written in terms of basic hypergeometric series as
\[  S_q(x) = - \frac{q}{1-q}\ {}_2\phi_1 \left( \genfrac{.}{.}{0pt}{}{q,q}{q^2}\ ;q,q \right)
             - \frac{qx}{(1-q)^2} \ {}_2\phi_2 \left( \genfrac{.}{.}{0pt}{}{q,q} {q^2,q^2}\ ; q,q^2x \right). \]
The growth of the coefficients in this Taylor series again shows that $S_q$ is an entire function of
order zero if we use the formula (\ref{eq:boas}) for the order of $\sum_{n=0}^\infty a_n z^n$, see
also \cite[\S 4]{gautschi}.

Next we mention the following $q$-integral representation, where we use 
Jackson's $q$-integral, see \cite[\S 1.11]{GR}
\begin{equation}\label{eq:defJacksonintegral}
\int_0^a f(t)\, d_qt = (1-q)a\, \sum_{k=0}^\infty f(aq^k)\, q^k,
\end{equation}
defined for functions $f$ whenever the right hand side converges.

\begin{property}\label{prop:qintegralrepr} For every $x \in \mathbb{C}$ we have
\[  S_q(x) = -\frac{q(1-x)}{1-q}\int_0^1 G_q(qx,qt)\, d_qt , \] 
with
\[ G_q(x,t) = \sum_{k=0}^\infty t^k (x;q)_k  
  = {}_2\phi_1 \left( \genfrac{.}{.}{0pt}{}{x,q}{0}\ ;q,t \right)
  = \frac{1}{1-t}\ {}_1\phi_1\left( \genfrac{.}{.}{0pt}{}{q}{qt}\ ; q,xt \right). \]
\end{property}

Since $\int_0^a f(t)\, d_qt\to \int_0^a f(t)\, dt$ when $q \to 1$ and
$G_q(x,t)\to 1/(1-t(1-x))$ when $q \to 1$ for $x > 0$, we see (at least formally) 
that Property \ref{prop:qintegralrepr} 
is a $q$-analogue of the integral representation
\[   \log(x) = - \int_0^1 \frac{1-x}{1-t(1-x)} \, dt, \qquad  x \notin (-\infty,0] \]
for the logarithm. 

\begin{proof}
Observe that 
\[ \frac{1-q}{1-q^{k+1}}= (1-q) \sum_{p=0}^\infty q^{(k+1)p} 
=\int_0^1 t^k\, d_qt. \] 
Inserting this in the definition 
\eqref{eq:S} of $S_q$ and interchanging summations, which is justified by the absolute 
convergence of the double sum,  gives the result. 
The identity between the basic hypergeometric series representing $G_q(x,t)$
is the case $c=0$ of \cite[(III.4)]{GR}.
\end{proof}

Note that, as in the proof of Property \ref{prop:TaylorS}, 
one can show that 
\begin{equation}\label{eq:altdefG}
G_q(x,t) = \sum_{j=0}^\infty \frac{(-xt)^j q^{j(j-1)/2}}{(t;q)_{j+1}}.
\end{equation}

\section{$q$-difference equation}\label{sec:qdiffeq}

The function $S_q$ satisfies a simple $q$-difference equation:

\begin{property}
The $q$-logarithm (\ref{eq:S}) satisfies
\begin{equation}  \label{eq:qrecur}
    S_q(x/q) - S_q(x) = 1 - (x;q)_\infty. 
\end{equation}
\end{property}
\begin{proof}
Recall the $q$-difference operator
\[    D_qf(x) = \frac{f(qx)-f(x)}{x(q-1)}, \]
then a simple exercise is 
\[    D_{1/q}(x;q)_k = - \frac{1-q^k}{1-q} (x;q)_{k-1}.  \]
Use this in (\ref{eq:S}) to find
\[   D_{1/q}  S_q(x) = \sum_{k=1}^\infty \frac{q^k}{1-q^k} \frac{1-q^k}{1-q} (x;q)_{k-1} 
    = \frac{q}{1-q} \sum_{k=0}^\infty q^k (x;q)_k.  \]
Observe that
$(x;q)_{k+1} - (x;q)_k = (x;q)_k [ 1-x q^k -1] = -xq^k (x;q)_k$, 
and summing  we find
$-x \sum_{k=0}^n q^k (x;q)_k = (x;q)_{n+1} -(x;q)_0$, 
and when $n \to \infty$
\[   \sum_{k=0}^\infty q^k (x;q)_k = \frac{1-(x;q)_\infty}{x} .  \]
If we use this result, then
\[    D_{1/q}  S_q(x) = \frac{q}{1-q} \frac{1-(x;q)_\infty}{x}, \]
which is \eqref{eq:qrecur}.
\end{proof}

In order to see how this is related to the classical derivative of $\log x$, one may
rewrite this as
\[    D_q((1-q) S_q(x)) = \frac{1}{x} - \frac{(qx;q)_\infty}{x}. \]
This $q$-difference equation can already be found in \cite[\S 6]{euler}, where Euler writes
$s=S_q(x)$ and $t=S_q(x/q)$ and gives the relation
\[  1+s-t=(1-x)\left(1-\frac{x}{a}\right)\left(1-\frac{x}{a^2}\right)\left(1-\frac{x}{a^3}\right)
    \left(1-\frac{x}{a^4}\right)\left(1-\frac{x}{a^5}\right) \cdots, \]
where $q=1/a$. 

As a corollary one has \cite[\S 7]{euler}

\begin{property}\label{prop:Sq-n}
For every positive integer $n$ one has
$S_q(q^{-n}) = n$.  
\end{property}

\begin{proof}
Use (\ref{eq:qrecur}) with $x=q^{-n+1}$ to find
$S_q(q^{-n}) - S_q(q^{-n+1}) = 1$,  
since $(x;q)_\infty$ vanishes whenever $x=q^{-n}$ for $n \geq 0$. The result then follows by induction
and $S_q(1) = 0$.
\end{proof}
It is this property, which is quite similar to $\log_a a^n = n$, where $\log_a$ is the logarithm with base $a$,
which gives $S_q$ the flavor of a $q$-logarithm, and which made
Euler consider this function as an interpolation of the logarithm, see
\cite[\S 1]{gautschi}. Observe that this interpolation property can be stated as follows: $-\log q\ S_q(x)$
approximates $\log x$ as $q \uparrow 1$ and for fixed $q$ this approximation is perfect 
if $x=q^{-n}$ $(n=1,2,\ldots)$.

Another interesting value is 
\[    S_q(0) = -\sum_{k=1}^\infty \frac{q^k}{1-q^k} = - \zeta_q(1), \]
which is a $q$-analogue of the harmonic series, where the $q$-analogue of the
$\zeta$-function is defined by
\[   \zeta_q(s) = \sum_{n=1}^\infty \frac{n^{s-1}q^n}{1-q^n}. \]
It has been proved, see
Erd\H{o}s \cite{erdos}, Borwein \cite{borwein1,borwein2}, 
Van Assche \cite{wva} 
that this quantity is irrational whenever $q=1/p$ with $p$ an integer $\geq 2$.
For the specific argument $1$ this coincides, up to a factor, 
with the value at $1$ of the $q$-$\zeta$-function considered 
by Ueno and Nishizawa \cite{UenoN}. 

The values of $S_q(q^n)$ for $n \in \mathbb{N}$ are distinctly different and for these
values we do not get the same flavor as the logarithm.

\begin{property}
For every positive integer $n$ one has
\begin{equation} \label{eq:Sqn}
    S_q(q^n) = -n + (q;q)_\infty \sum_{k=0}^{n-1} \frac{1}{(q;q)_k}. 
\end{equation}
\end{property}
\begin{proof}
Choose $x=q^{k+1}$ in (\ref{eq:qrecur}), then
$S_q(q^k) - S_q(q^{k+1}) = 1-(q^{k+1};q)_\infty$. 
Summing and the telescoping property gives
\[  S_q(q^0) - S_q(q^n) = \sum_{k=0}^{n-1} \left( S_q(q^k) - S_q(q^{k+1}) \right) = 
 n - \sum_{k=0}^{n-1} (q^{k+1};q)_\infty. \]
By Property \ref{prop:Sq-n}  we have $S_q(1)=0$. Now
$(q^{k+1};q)_\infty = (q;q)_\infty/(q;q)_k$ gives the required expression \eqref{eq:Sqn}.
\end{proof}

In order to see how this approximates $\log x$, one may reformulate this as
\[    -\log q \ S_q(q^n) = \log q^n - \log q \sum_{k=0}^{n-1} (q^{k+1};q)_\infty . \]
In \cite[\S 10]{euler} Euler writes $s=S_q(q^n)$, $t=S_q(q^{n-1})$, $u=S_q(q^{n-2})$ and he writes
the recursion
\[     s= \frac{2t-u+aq^n(1-t)}{1-aq^n}, \]
where $q=1/a$. In contemporary notation we write $y_n = S_q(q^n)$ and obtain the recurrence
relation
\[     y_n (1-q^{n-1}) - (2-q^{n-1})y_{n-1} +y_{n-2} = q^{n-1}.  \]
One can verify that this recurrence relation indeed holds for $y_n=S_q(q^n)$ given in (\ref{eq:Sqn}). 
More generally one in fact has
\[  (1-qx)S_q(q^2x) - (2-qx)S_q(qx) + S_q(x) = qx, \]
which is non-homogeneous a second order $q$-difference equation for $S_q$.

Note that the explicit evaluation $S(q^{-n})=n$, $n\in\N$, gives the
following summation formulas
\begin{equation}\label{eq:sumformfromS1}
\sum_{k=1}^n \frac{(q^{-n};q)_k}{1-q^k} q^k =-n, 
\qquad
\sum_{k=1}^\infty \frac{q^{k(k+1)/2} (-1)^{k-1} q^{-nk}}{(1-q^k)\, (q;q)_k} 
= n + \sum_{k=1}^\infty \frac{q^k}{1-q^k}
\end{equation}
using the definition of $S_q(x)$ and the Taylor expansion in 
Property \ref{prop:TaylorS}. Similarly, the evaluation at $q^n$, $n\in\N$, 
given in \eqref{eq:Sqn} gives the summation formulas
\begin{equation}\label{eq:sumformfromS2}
 \begin{split}
&\sum_{k=1}^\infty \frac{(q^n;q)_k}{1-q^k} q^k = n - \sum_{k=0}^{n-1} (q^{k+1};q)_\infty, \\
&\sum_{k=1}^\infty \frac{q^{k(k+1)/2} (-1)^{k-1} q^{nk}}{(1-q^k)\, (q;q)_k} 
= - n + \sum_{k=1}^\infty \frac{q^k}{1-q^k} + \sum_{k=0}^{n-1} (q^{k+1};q)_\infty.
\end{split}
\end{equation}
Note that all infinite series are absolutely convergent 
and that for $n=0$ the results in \eqref{eq:sumformfromS1} and
\eqref{eq:sumformfromS2} coincide. The first sums become trivial, and the
second gives an expansion for the $\zeta_q(1)$
\begin{equation}\label{eq:altexpreszeta1}
\zeta_q(1) = \sum_{k=1}^\infty \frac{q^k}{1-q^k} = 
\sum_{k=1}^\infty \frac{q^{k(k+1)/2} (-1)^{k-1}}{(1-q^k)\, (q;q)_k}. 
\end{equation}
Using \eqref{eq:altexpreszeta1} in Property \ref{prop:TaylorS} gives the
 expansion
\[
 S_q(x) = 
-\sum_{k=1}^\infty \frac{q^{k(k+1)/2} (-1)^{k-1}(1-x^k)}{(1-q^k)\, (q;q)_k},
\]
so that in particular
\[
- \frac{dS_q}{dx}(1) = \lim_{x\to 1} \frac{S_q(x)}{1-x} 
= -\sum_{k=1}^\infty \frac{k\, q^{k(k+1)/2} (-1)^{k-1}}{(1-q^k)\, (q;q)_k}.
\]

\section{An extension of the $q$-logarithm and Lambert series}\label{sec:extensionqlog}

Having the definition of $S_q(x)$ resembling Lambert series, it is natural
to look for the extension
\begin{equation}\label{eq:defF}
F_q(x,t) = - \sum_{k=1}^\infty (x;q)_k \frac{t^k}{1-t^k},
\end{equation}
which is a Lambert series, see \cite[\S 58.C]{knopp}. 
Since $|(x;q)_k|\leq (-|x|;|q|)_k\leq (-r;|q|)_\infty$
for $x$ in $\{x\in\C\mid |x|\leq r\}$, the convergence in 
\eqref{eq:defF} is uniform on compact sets in $x$ and on compact subsets of the open
unit disk in $t$. Also since the series 
$-\sum_{k=1}^\infty (x;q)_k t^k$ is absolutely convergent for $|t|<1$ 
uniformly in $x$ in compact sets, it follows
by \cite[Satz 259]{knopp}, that $F_q$ is analytic for $(x,t)\in \C\times \{t\in\C\mid |t|<1\}$.
Observe that $S_q(x)= F_q(x,q)$. 

The general theory of Lambert series then gives the power series 
of $F$ in powers of $t$;
\[
F_q(x,t) = \sum_{\ell=1}^\infty \Bigl( \sum_{k | \ell} (x;q)_k\Bigr)  t^\ell 
\ \Longrightarrow\  S_q(x) =\sum_{\ell=1}^\infty \Bigl( \sum_{k | \ell} (x;q)_k\Bigr)  q^\ell
\]
We are mainly interested in 
the power series development with respect to $x$. 

\begin{property}\label{prop:expansionFinpowersx}
For $|t|<1$ one has
\[
F_q(x,t) = - \sum_{k=1}^\infty \frac{t^k}{1-t^k} 
- \sum_{\ell=1}^\infty x^\ell (-1)^\ell q^{\ell(\ell-1)/2} \left( \sum_{n=1}^\infty 
t^{n\ell}\frac{(t^n q^{\ell+1};q)_\infty}{(t^n;q)_\infty}\right). 
\]
\end{property}

In case $t=q$, Property \ref{prop:expansionFinpowersx} 
reduces to Property \ref{prop:TaylorS}, and this is equivalent to the summation
formula
\begin{equation}\label{eq:nontrivsummation}
\sum_{n=1}^\infty q^{n\ell} \frac{(q^{\ell+n+1};q)_\infty}{(q^n;q)_\infty}  
= \frac{q^\ell}{(1-q^\ell)\, (q;q)_\ell}\ 
\Longrightarrow\ 
\sum_{n=1}^\infty \frac{(q;q)_{n-1}}{(q^{\ell+1};q)_n} q^{n\ell} = \frac{q^\ell}{1-q^\ell}
\end{equation}
for $\ell\in\N$, $\ell\geq 1$. This can be obtained as a 
special case of $q$-Gauss sum
\cite[(1.5.1)]{GR}.

\begin{proof} The proof is along the same lines as the proof
of Property \ref{prop:TaylorS}. We find similarly 
\begin{equation*}
F_q(x,t) = -\sum_{k=1} \frac{t^k}{1-t^k} 
- \sum_{j=1}^\infty q^{j(j-1)/2} (-xt)^j \sum_{\ell=0}^\infty 
\frac{(q^{j+1};q)_\ell}{(q;q)_\ell} \frac{t^\ell}{1-t^{j+\ell}}
\end{equation*}
and we write
\begin{multline*}
\sum_{\ell=0}^\infty 
\frac{(q^{j+1};q)_\ell}{(q;q)_\ell} \frac{t^\ell}{1-t^{j+\ell}} = 
 \sum_{\ell=0}^\infty 
\frac{(q^{j+1};q)_\ell}{(q;q)_\ell} t^\ell \sum_{p=0}^\infty t^{p(j+\ell)} 
\\ = \sum_{p=0}^\infty t^{jp} \sum_{\ell=0}^\infty 
\frac{(q^{j+1};q)_\ell}{(q;q)_\ell} t^{\ell(1+p)}
= \sum_{p=0}^\infty t^{jp} \frac{(t^{1+p}q^{j+1};q)_\infty}{(t^{1+p};q)_\infty}
\end{multline*}
using the $q$-binomial theorem again and the absolute convergence of the double sum,
which justifies the interchange of summations. Using this and replacing $n=p+1$
gives the result. 
\end{proof}

Consider the
case $t=q^2$. Following the line of proof of Property 
\ref{prop:TaylorS} we write
\[
-\sum_{k=1}^\infty \frac{q^{2k}(x;q)_k}{1-q^{2k}} 
= - \sum_{k=1}^\infty \frac{q^{2k}}{1-q^{2k}} - 
\sum_{j=1}^\infty \frac{(-1)^j q^{j(j-1)/2} x^j}{(q;q)_j} 
\sum_{\ell=0}^\infty \frac{(q;q)_{\ell+j} \, q^{2\ell+2j}}{(q;q)_\ell \, (1-q^{2\ell+2j})}
\]
and the inner sum over $\ell$ can be written as
\[
\sum_{\ell=0}^\infty \frac{(q;q)_{\ell+j-1}\, q^{2\ell+2j}}{(q;q)_\ell\, (1+q^{\ell+j})} =
\frac{(q;q)_{j-1}q^{2j}}{1+q^j} \sum_{\ell=0}^\infty
\frac{(q^j;q)_\ell (-q^j;q)_\ell}{(q;q)_\ell (-q^{j+1};q)_\ell} q^{2\ell}.
\] 
Using Property \ref{prop:expansionFinpowersx} for $t=q^2$ then gives
\begin{equation}\label{eq:resultz=q2}
\sum_{n=1}^\infty q^{2nj} \frac{(q^{2n+j+1};q)_\infty}{(q^{2n};q)_\infty}
= \frac{q^{2j}}{(1-q^{2j})} 
\sum_{\ell=0}^\infty
\frac{(q^j;q)_\ell (-q^j;q)_\ell}{(q;q)_\ell (-q^{j+1};q)_\ell} q^{2\ell}. 
\end{equation}
This can also be proved directly using the $q$-binomial theorem and
geometric series. 
We can rewrite \eqref{eq:resultz=q2} in standard basic hypergeometric
series form, see \cite{GR}, as the quadratic transformation 
\begin{equation}\label{eq:resultz=q2-v2}
\frac{(1-q^{2j})}{(q^2;q)_{j+1}}\,_3\phi_2 \left( \genfrac{.}{.}{0pt}{}{q^2,q^2,q^3}
{q^{j+3}, q^{j+4}}\ ;q^2, q^2 \right)
=  \,_2\phi_1 \left( \genfrac{.}{.}{0pt}{}{q^j,-q^j}
{-q^{j+1}}\ ;q, q^2 \right). 
\end{equation}

Analogous to Property \ref{prop:qintegralrepr}, and using the notation
of Property \ref{prop:qintegralrepr} we have the following.

\begin{property}\label{prop:qintegralFrepr}
For $|p|<1$ one has 
\[ F_q(x,p) = \frac{-p(1-x)}{(1-p)} \int_0^1 G(qx,pt)\, d_pt. \]
\end{property}


\section{A $q$-analogue of the dilogarithm}\label{sec:qdilog}

Euler's dilogarithm is defined by the first equality in 
\[
\Li (x)  = \sum_{n=1}^\infty \frac{x^n}{n^2} = - \int_0^x \frac{\log(1-t)}{t}\, dt 
= - \int_{1-x}^1 \frac{\log(t)}{1-t}\, dt = \frac{\pi^2}{6} - \Li(1-x)
\]
for $0\leq x\leq 1$, 
see \cite{lewin}, \cite{kirillov}, for more information and references. 
Here we use $\Li (1) = \zeta(2) = \frac{\pi^2}{6}$. 
In particular, $x \frac{d\Li}{dx} = -\log(1-x)$, and the definition
by the series can be extended to complex $x$ being absolutely convergent
for $|x|\leq 1$. 
 
We define the $q$-dilogarithm by
\begin{equation}\label{eq:defqdilog}
\Li (x;q) = \sum_{k=1}^\infty \frac{q^k}{(1-q^k)^2}(x;q)_k.
\end{equation}
We have $\lim_{q\uparrow 1} (1-q)^2 \Li (x;q) = 
\sum_{k=1}^\infty (1-x)^k/k^2 = \Li (1-x)$. 
In this case we can justify the interchange of the limit and
summation using dominated convergence. We assume $0<q<1$, and 
we first observe that
$|(x;q)_k|\leq 1$ for $|1-x|\leq 1$. Next we use
\[   \frac{1-q^k}{1-q} = \sum_{j=0}^{k-1} q^j = 
  q^{(k-1)/2} \begin{cases}
       \sum_{j=0}^{\frac{k}{2}-1} \left( q^{j+\frac12}+q^{-j-\frac12} \right), &  \textrm{$k$ even}, \\
       1 + \sum_{j=0}^{\frac{k-1}{2}-1} \left(  q^{j+1} + q^{-j-1} \right), & \textrm{$k$ odd}, 
       \end{cases}   \]
and $x+1/x \geq 2$ for $x \in [0,1]$ then gives 
\[ \frac{1-q^k}{1-q} \geq k q^{(k-1)/2}, \]
so that
\[   q^k \frac{(1-q)^2}{(1-q^k)^2} \leq \frac{1}{k^2}. \]
Combining both estimates gives 
\[ |\frac{q^k}{(1-q^k)^2}(x;q)_k| \leq \frac{1}{k^2} \]
 for $|1-x|\leq 1$ and dominated convergence is established.
%

We list some properties of the $q$-dilogarithm.
In the following we use $\zeta_q(2)=\sum_{k=1}^\infty \frac{q^k}{(1-q^k)^2}$,
as an analogue of $\frac16 \pi^2$. This is equal to the $q$-$\zeta$-function
\[ \zeta_q(s) = \sum_{n=1}^\infty \frac{n^{s-1}q^n}{1-q^n} \] 
for $s=2$ since
\[  \sum_{n=1}^\infty \frac{nq^n}{1-q^n} = \sum_{n=1}^\infty \sum_{k=1}^\infty nq^{nk}
    = \sum_{k=1}^\infty \frac{q^k}{(1-q^k)^2}, \]
(see, e.g., \cite[Part VIII, Chapter 1, problem 75]{polszego}). This quantity
was considered by Zudilin \cite{zud,zud2}, Krattenthaler et al. \cite{krat}, 
Postelmans and Van Assche \cite{pos}, who studied its
irrationality when $1/q$ is an integer $\geq 2$. Note that this does no longer
correspond to Ueno and Nishizawa \cite{UenoN}, who essentially have 
$\sum_{k=1}^\infty \frac{q^{2k}}{(1-q^k)^2}$ as the value at $2$ for their
$q$-$\zeta$-function.

\begin{property}\label{prop:Li} $\Li (\,\cdot\,;q)$ is an entire function of
order zero. Moreover, we have the special values 
\[
\Li (1;q)=0, \quad \Li (0;q) =  \zeta_q(2), 
\quad \Li (q^{-n};q) = -\sum_{k=1}^n \frac{k}{1-q^k},
\] 
and
$(1-q)(1-x)\, \bigl( D_q \Li (\,\cdot\,;q)\bigr) (x) = S_q(x)$ and
\[
\Li (x;q) = \zeta_q(2) + \frac{1}{1-q} \int_0^x \frac{S_q(t)}{1-t}\, d_q t.
\]
Moreover, the $q$-dilogarithm has the Taylor expansion
\[
\Li (x;q) = \zeta_q(2) +  \sum_{j=1}^\infty \frac{(-1)^j q^{j(j+1)/2}\, x^j}{(1-q^j)^2}
\,_2\phi_1 \left( \genfrac{.}{.}{0pt}{}{q^j,q^j}{q^{j+1}}
\ ;q, q \right).
\]
\end{property}

Here the ${}_2\phi_1$-series is defined by
\[
\,_2\phi_1 \left( \genfrac{.}{.}{0pt}{}{q^j,q^j}{q^{j+1}}
\ ;q, q \right) = \sum_{\ell=0}^\infty 
\frac{(q^j;q)_\ell (q^j;q)_\ell}{(q;q)_\ell (q^{j+1};q)_\ell} \, q^\ell 
= \sum_{\ell=0}^\infty 
\frac{(q^j;q)_\ell (1-q^j)}{(q;q)_\ell (1-q^{j+\ell})} \, q^\ell .
\]
Unfortunately, this series cannot be summed using the (non-terminating) 
$q$-Chu-Vandermonde sum.

Note that after multiplying the integral representation for $\Li (x;q)$ by $(1-q)^2$
we can take  a formal limit $q\uparrow 1$ to get
\[
\Li (1-x) = \frac{\pi^2}{6} + \int_0^x \frac{\log(t)}{1-t} \, dt 
= - \int_0^{1-x} \frac{\log(1-t)}{t}\, dt,
\]
so that we recover the integral representation for the
dilogarithm. 

\begin{proof}
The proof of $\Li(\,\cdot\,;q)$ being an entire function of order zero
 is derived as in Property \ref{prop:entireS}.
Since
$(qx;q)_k-(x;q)_k= x(1-q^k)(qx;q)_{k-1}$ we obtain
\begin{equation}\label{eq:qdiffdilog}
\begin{split}
\Li (qx;q)-\Li (x;q) = \frac{x}{1-x} \sum_{k=1}^\infty 
\frac{q^k (x;q)_k}{1-q^k} = \frac{-x}{1-x} \, S_q(x).
\end{split}
\end{equation}
This implies $(1-q)(1-x)\, \bigl( D_q \Li (\,\cdot\,;q)\bigr) (x) = S_q(x)$.

Using \eqref{eq:qdiffdilog} for $x=q^{-n}$, $n\in\N$, and $\Li (1;q)=0$, 
$S(q^{-n})=n$ we find
the value for 
$\Li (q^{-n};q)$.  
Iterating \eqref{eq:qdiffdilog} we get
\[
\Li (x;q) = \sum_{k=0}^N \frac{xq^k}{1-xq^k} S_q(xq^k) + \Li (xq^{N+1};q)
\]
and by letting $N\to\infty$ we get the convergent series expansion
\[
\Li (x;q) =  \Li (0;q) + \sum_{k=0}^\infty \frac{xq^k}{1-xq^k} S_q(xq^k) 
= \zeta_q(2) + \frac{1}{1-q} \int_0^x \frac{S_q(t)}{1-t} \, d_qt.
\]

Finally, the Taylor expansion proceeds as in the proof of
Property \ref{prop:TaylorS}, and we find
\begin{equation*}
\Li (x;q) = \sum_{k=1}^\infty \frac{q^k}{(1-q^k)^2} 
+ \sum_{j=1}^\infty \frac{(-x)^j q^{j(j-1)/2}}{(q;q)_j}
\sum_{\ell=0}^\infty \frac{(q;q)_{j+\ell} \, q^{j+\ell}}{(q;q)_\ell\, (1-q^{j+\ell})^2}
\end{equation*}
The inner sum over $\ell$ can be rewritten as
\begin{equation*}
\frac{q^j (q;q)_{j-1}}{1-q^j} \sum_{\ell=0}^\infty 
\frac{(q^j;q)_\ell (q^j;q)_\ell}{(q;q)_\ell (q^{j+1};q)_\ell} \, q^\ell 
\end{equation*}
and this gives the result.
\end{proof}

The evaluation of the $q$-dilogarithm gives the following
summation, cf. \eqref{eq:sumformfromS1},
\begin{multline}\label{eq:sumformfromLi}
\sum_{k=1}^n \frac{(q^{-n};q)_k q^k}{(1-q^k)^2} 
= - \sum_{k=1}^n \frac{k}{1-q^k} \\ = 
\sum_{j=1}^\infty \frac{q^j}{(1-q^j)^2} 
+ \sum_{j=1}^\infty \frac{(-1)^j q^{j(j+1)/2}\, q^{-nj}}{(1-q^j)^2}
\,_2\phi_1 \left( \genfrac{.}{.}{0pt}{}{q^j,q^j}{q^{j+1}}
\ ;q, q \right).
\end{multline} 
In particular, for $n=0$ we obtain an alternating series representation
for $\zeta_q(2)$;
\[
\zeta_q(2) = \sum_{j=1}^\infty \frac{(-1)^{j-1} q^{j(j+1)/2}}{(1-q^j)^2}
\,_2\phi_1 \left( \genfrac{.}{.}{0pt}{}{q^j,q^j}{q^{j+1}}
\ ;q, q \right).
\]

Writing $\Li (x;q) = \sum_{n=0}^\infty a_n x^n$, $S_q(x) = \sum_{n=0}^\infty b_n x^n$ 
temporarily, then \eqref{eq:qdiffdilog} implies that 
$q^na_n-a_n$ equals the coefficient, say $c_n$, of $x^n$ in 
$- S_q(x) x/(1-x)$. Using $-x/(1-x) = \sum_{k=1}^\infty -x^k$, it follows
that $c_n= -\sum_{p=0}^{n-1} b_p$. Note that the relation is trivial in case
$n=0$, and for integer $n\geq 1$ we find from the explicit Taylor expansions
for $S_q(\,\cdot\,)$ and $\Li (\,\cdot\,;q)$ the relation
\begin{equation*}
\frac{(-1)^{n-1} q^{n(n+1)/2}}{(1-q^n)}
\,_2\phi_1 \left( \genfrac{.}{.}{0pt}{}{q^n,q^n}{q^{n+1}}
\ ;q, q \right)
=  \sum_{k=1}^\infty \frac{q^k}{1-q^k} + \sum_{j=1}^{n-1} 
 \frac{(-1)^jq^{j(j+1)/2}} {(1-q^j)\, (q;q)_j}.
\end{equation*}
Note that this relation gives an explicit expression for the 
remainder if approximating $\zeta_q(1)$ with the alternating series
as in \eqref{eq:altexpreszeta1}. Of course, we get the same result if
we use the Taylor expansion of $S_q$ as in Property 
\ref{prop:TaylorS} in the integral representation for the
$q$-dilogarithm in Property \ref{prop:Li}. 

The classical dilogarithm satisfies many interesting 
properties, such as a simple functional equation, a
five-term recursion, a characterisation by these 
first two properties, explicit evaluation at certain 
special points, etc., see \cite{lewin}, \cite{kirillov}
for more information and references. 
It would be interesting to see if these interesting
properties have appropriate analogues for the
$q$-analogue of the dilogarithm discussed here.

\section{Other $q$-logarithms}\label{sec:otherqlogs}

In physics literature, see e.g. \cite{tsallis}, one defines 
$\ln_q(x) = \frac{x^{1-q}-1}{1-q}$.  
There are no $q$-series, $q$-Pochhammer symbols, $q$-difference
relations, etc. The choice of the letter $q$ and the fact that
$\lim_{q \to 1} \ln_q(x)  = \log x$
is not sufficient motivation to call this a $q$-analogue. It just shows that the logarithmic function
is somewhere between the constant function and powers $x^\alpha-1$ for $\alpha > 0$.

Borwein \cite{borwein2}, Zudilin et al. \cite{matvazud}, Van Assche \cite{wva} 
consider 
\[   \ln_q(1+z) = \sum_{k=1}^\infty \frac{(-1)^{k} z^k}{1-q^k}, \qquad |z| < |q|, \]
with $|q| > 1$. They prove that $\ln_q(1+z)$ is irrational for $z=\pm1$ and $q$ an integer greater than $2$.
For $z=-1$ one has a $q$-analogue of the harmonic series and this is essentially
the generating function of $d_n = \sum_{k | n} 1$, i.e. the number of divisors of $n$.
A similar formula, but now for $0<q<1$ 
\[
\log_q(z) = \sum_{k=1}^\infty \frac{z^n}{1-q^n} = \frac{z \, e_q'(z)}{e_q(z)}, \qquad |z| <1
\]
has been considered as a $q$-analogue of the logarithm by 
Kirillov \cite{kirillov} and Koornwinder \cite{koornwinder}. 
This $q$-analogue is well adapted to non-commutative algebras,
see \cite[\S 2.5, Ex.~11]{kirillov}, \cite[Prop.~6.1]{koornwinder},
since $\log_q(x+y-xy)=\log_q(x)+\log_q(y)$ for $xy=qyx$. 
The corresponding $q$-analogue of the dilogarithm,
provisionally denoted by $\widetilde\Li (x;q)$, is defined
by 
\[
\widetilde\Li (x;q) = \sum_{k=1}^\infty \frac{z^k}{k\, (1-q^k)} = 
\log (e_q(z))\  \Longrightarrow\ \log_q(z) = z\, \widetilde\Li' (z;q). 
\]
Zudilin \cite{zud2} considers a similar $q$-logarithm but a different $q$-dilogarithm
\[  L_1(x;q) = \sum_{n=1}^\infty \frac{(xq)^n}{1-q^n}, \quad 
    L_2(x;q) = \sum_{n=1}^\infty \frac{n(xq)^n}{1-q^n}, \]
and mainly studies simultaneous rational approximation to $L_1$ and $L_2$
in order to obtain quantitative linear independence over $\mathbb{Q}$
for certain values of these functions.

Other $q$-logarithms are defined as inverses of $q$-exponential functions, see
Nelson and Gartley \cite{nelson} for two different cases viewed from 
complex function theory, and Chung et al. \cite{chung}, where implicitly $q$-commuting variables are used.
Fock and Goncharov \cite{FockG,gon} introduce a $q$-logarithm of $\ln(e^z+1)$ by an integral.
The corresponding $q$-dilogarithm is essentially Ruijsenaars' hyperbolic
$\Gamma$-function, see \cite[II.A]{Ruis}. 
For other $q$-logarithms based on Jacobi theta functions, see Sauloy 
\cite{sauloy} and Duval \cite{duval}, where the $q$-logarithms play a role
in difference Galois theory in constructing the analogue of a unipotent 
monodromy representation.

\bibliographystyle{amsplain}

\end{document}